\documentclass[11pt]{article}
\usepackage{epsfig}
\usepackage{graphics}
\usepackage{graphicx}
\usepackage{float}
\usepackage{amsthm}
\usepackage{amsmath,amssymb,amsfonts}
\newtheorem{theorem}{Theorem}[section]
\newtheorem{definition}{Definition}[section]
\newtheorem{lemma}{Lemma}[section]

\newtheorem{proposition}{Proposition}[section]

\newcommand{\newsection}[1]{\section{#1}\setcounter{equation}{0}}

\newcommand{\vs}{\vspace{4.6mm}}

\def\b#1{{\bf #1}}

\def \defi{\stackrel{\scriptstyle{{\rm def}}}{=}}

\begin{document}

\begin{center} {\rm \bf \LARGE On non-trivial barrier solutions of the dividend problem for a diffusion
under constant and proportional transaction costs}\end{center} \vspace{0.1cm}
\begin{center}Lihua Bai$^{\mbox{a}}$ and Jostein 
Paulsen$^{\mbox{b}}$\footnote{Corresponding author} \end{center}
\begin{center}
$^{\mbox{a}}$Department of Mathematical Sciences and LPMC, Nankai                                                                                                                                                                                                                                                                       
University, China\\
$^{\mbox{b}}$Department of Mathematical Sciences, University of Copenhagen, Denmark \end{center}

\vspace{0.2cm}
\begin{abstract}
In Bai and Paulsen (SIAM J. Control optim. 48, 2010) the optimal dividend problem under  transaction costs was analyzed for a rather general class
of diffusion processes. It was divided  into several subclasses, and for the majority of subclasses the optimal
policy is a simple barrier policy; whenever the process hits an upper barrier
$\bar{u}^*$, reduce it to $\bar{u}^*-\xi$ through a dividend payment. 
After transaction costs, the shareholder receives $k\xi-K$. 

It was  proved that  a simple barrier strategy is not always optimal, and
here  these more difficult cases are solved. The optimal solutions are rather complicated, but  interesting.

\end{abstract}


\vspace{0.2cm}

\newsection{Introduction and problem formulation}
Let
 $(\Omega, {\cal F}, \{{\cal F}_{t}\}_{t\geq 0},
 P)$ be a probability space satisfying the usual conditions, i.e. the
 filtration $\{{\cal F}_{t}\}_{t\geq 0}$ is right continuous and
 $P$-complete. Assume that the uncontrolled surplus process follows
 the stochastic differential equation
 \begin{eqnarray*}
dX_{t} =\mu(X_{t})dt+\sigma(X_{t})dW_{t},\quad X_{0}=x,
\end{eqnarray*}
where $W$ is a Brownian motion and $\mu(x)$
and $\sigma(x)$ are Lipschitz-continuous. Let the  company
pay dividends to its shareholders, but at a fixed transaction cost
$K>0$ and a tax rate $1-k<1$, so that $k>0$. We will allow $k>1$, opening up for other
interpretations than that $1-k$ is a tax rate. This means that if $\xi>0$ is the
amount the capital is reduced by, then the net amount of money  the
shareholders receive is $k\xi-K$. It can be argued that taxes are paid on dividends after 
costs, so an alternative would be to use $k(\xi-K)=k\xi-kK$, but clearly 
this is just a reparametrization. Furthermore, different investors may 
have different tax rates, so $1-k$ should be interpreted as an average tax rate.

Since every dividend payment results in a fixed
transaction cost, the  company should not pay out
dividends continuously, but only at discrete time epochs. Therefore,
a strategy can be described by
$$\pi=(\tau_{1}^{\pi},\tau_{2}^{\pi},.\ .\
., \tau_{n}^{\pi},.\ .\ .; \xi_{1}^{\pi},\xi_{2}^{\pi},.\ .\ .,\xi_{n}^{\pi},.\ .\ .),$$
 where  $\tau_{n}^{\pi}$  and $\xi_{n}^{\pi}$ denote the times and
amounts of dividends. Thus, when applying the strategy $\pi$,
   the resulting surplus process ${X_{t}^{\pi}}$ is
   given by
 \begin{eqnarray*}
 X_{t}^{\pi}=x+\int_{0}^{t}\mu(X_{s}^{\pi})ds+\int_{0}^{t}\sigma(X_{s}^{\pi})dW_{s}-\mathop
{\sum}_{n=1}^{\infty}1_{\{\tau_{n}^{\pi}<t\}}\xi_{n}^{\pi}.
\end{eqnarray*}
Note that this makes $X^{\pi}$ left continuous, so that $\xi_n^{\pi}=X_{\tau_n^{\pi}}^{\pi}
-X_{\tau_n^{\pi}+}^{\pi}$.

\begin{definition}
A strategy $\pi$ is said to be
admissible  if
\begin{itemize}
\item[(i)] $ 0\leq \tau_{1}^{\pi}$ and for $n\geq 1$, $\tau_{n+1}^{\pi}>\tau_n^{\pi}$ on
$\{\tau_n^{\pi}<\infty\}$.
\item[(ii)] $\tau_{n}^{\pi}$ is a stopping time with respect to  $\{{\cal
F}_{t}\}_{t\geq 0},$ $ n=1,2.\ .\ .\ .$
\item[(iii)] $\xi_{n}^{\pi} $ is measurable with respect to ${\cal
 F}_{\tau_{n}^{\pi}+},\ n=1,2.\ .\ .\ .$
\item[(iv)] $\tau_n^{\pi}\rightarrow\infty\;\;\mbox{a.s. as}\;\;n\rightarrow \infty$.
\item[(v)] $0<\xi_{n}^{\pi}\leq X_{\tau_{n}}^{\pi}$.
\end{itemize}
 We denote the set of all admissible strategies by $\Pi.$ 
\end{definition}

Another natural admissibility condition is that net money received should be positive,
that is $k\xi -K>0$.  However, we are looking for optimal policies, and a policy
that allows $k\xi-K\leq 0$ can never be optimal, so it can be dropped as a condition.

 With each admissible strategy $\pi$ we define the
corresponding
 ruin time as
 \begin{eqnarray*}
 \tau^{\pi}={\rm inf}\{t\geq 0:\ X_{t}^{\pi}<0\},
  \end{eqnarray*}
and the performance function $V_{\pi}(x)$ as
$$V_{\pi}(x)=E_x\left[\mathop
{\sum}_{n=1}^{\infty}e^{-\lambda \tau_{n}^{\pi}}(k\xi_{n}^{\pi}-K)
1_{\{\tau_{n}^{\pi}\leq \tau^{\pi}\}}\right],$$
where by $P_x$ we mean the probability measure conditioned on $X_{0}=x$.
$V_{\pi}(x)$  represents the expected total
discounted dividends received by the shareholders until ruin 
when the initial reserve is $x$. \\
\indent Define the optimal return function by
 \begin{eqnarray*}
 V^{*}(x)=\mathop
 {\sup}_{\pi\in \Pi}V_{\pi}(x),
 \end{eqnarray*}
 and the optimal strategy, if it exists, by  $\pi^{*}$ so that  $V_{\pi^{*}}(x) =
 V^{*}(x)$.

Another natural admissibility condition is that net money received should be positive,
that is $k\xi -K>0$.  However, we are looking for optimal policies, and a policy
that allows $k\xi-K\leq 0$ can never be optimal, so it can be dropped as a condition.

\begin{definition}
A (simple) lump sum dividend barrier strategy
$\pi=\pi_{\bar{u},\underline{u}}$ with  parameters $0\leq\underline{u}<\bar{u}$,
is given by:
\begin{itemize}
\item When $X_t^\pi<\bar{u}$, do nothing.
\item When $X_t^\pi\geq \bar{u}$, reduce $X_t^\pi$ to $\underline{u}$ through a
dividend payment.
\end{itemize}

With a lump sum dividend barrier strategy $\pi_{\bar{u},\underline{u}}$,
the corresponding value function will be denoted by $V_{\bar{u},\underline{u}}(x)$.

A two-level lump sum dividend barrier strategy
$\pi=\pi_{(\bar{u}_1,\underline{u}_1)(u_2^c,\bar{u}_2,\underline{u}_2)}$ with  parameters 
$0\leq \underline{u}_1<\bar{u}_1<u_2^c <\bar{u}_2$ and $0\leq \underline{u}_2<\bar{u}_2$
is given by:
\begin{itemize}
\item When $X_t^\pi<\bar{u}_1$, do nothing.
\item When $\bar{u}_1\leq X_t^\pi \leq u_2^c$, reduce $X_t^\pi$ to $\underline{u}_1$ through
a dividend payment.
\item When $u_2^c<X_t^\pi<\bar{u}_2$, do nothing.
\item When $X_t^\pi \geq \bar{u}_2$, reduce $X_t^\pi$ to $\underline{u}_2$ through a dividend 
payment.
\end{itemize}

With a two-level lump sum dividend barrier strategy 
$\pi_{(\bar{u}_1,\underline{u}_1)(u_2^c,\bar{u}_2,\underline{u}_2)}$, the corresponding value function
will be denoted by $V_{(\bar{u}_1,\underline{u}_1)(u_2^c,\bar{u}_2,\underline{u}_2)}(x)$.
\end{definition}

\vspace{0.2cm}
\indent We will work under the following set of assumptions:
\begin{itemize}
\item[A1.] $\mid \mu(x)\mid +\mid \sigma(x)\mid\leq C(1+x)$ for all $x\geq
0$ and some $C>0.$
\item[A2.] $\mu(x)$ and $\sigma(x)$ are  continuously differentiable
 and Lipschitz continuous and
the derivatives $\mu'(x)$ and $\sigma'(x)$ are Lipschitz continuous for
all $x\geq 0.$
\item[A3.] $\sigma^{2}(x)>0$ for all $x\geq0.$ 
\item[A4.] There exists a number $x_{\lambda}\in [0,\infty)$ so that 
$\mu'(x)> \lambda$ for all $x< x_{\lambda}$ and 
$\mu'(x)\leq \lambda$ for all $x\geq x_{\lambda}$. The number $\lambda$ is a discounting rate.
\end{itemize}
Note that under A3, any dividend payment that reduces the capital to zero will result in ruin.

\vspace{0.2cm}
For  $g\in C^{2}(0,\infty),$  define the operator $L$ by
$$Lg(x)=\frac{1}{2}\sigma^{2}(x)g''(x)+\mu(x)g'(x)-\lambda g(x).$$
It is well known, see e.g. \cite{P3}, that under the
assumptions A1-A3 any solution of  $Lg(x)=0$ is in $C^{2}(0,\infty).$  
Let $g_1$ and $g_{2}$ be two
independent solutions of $Lg(x)=0$, chosen so that
$g(x)=g_{1}(0)g_{2}(x)-g_{2}(0)g_{1}(x)$ has $g'(0)>0$. Such a
solution will be called a canonical solution. Then any
solution of $LV(x)=0$ with $V(0)=0$ and $V'(0)>0$ is of the form
$$V(x)=cg(x),\qquad c>0.$$
Under Assumptions A1-A4, it was proved in \cite{Bai10b} that there are two basic 
possibilities for the canonical solution $g$.
\begin{itemize}
\item[P.] There is an $x^{*}$ with $0\leq x^*\leq \infty$
so that $g$ is concave on $[0,x^{*}]$ and convex on
$(x^{*},\infty)$. In particular $x^{*}=0$ if and only if
$\mu(0)\leq 0$, and by the definition of $x^*$, $g''(x^{*})=0$ when
$0<x^{*}<\infty.$
\item[R.] There are numbers $x_1^*$ and $x_2^*$ with $0<x_1^*\leq x_{\lambda}\leq x_2^*\leq \infty$ so
that $g$ is convex on $[0,x_1^*)$, concave on $[x_1^*,x_2^*)$ and, if $x_2^*<\infty$,
convex on $[x_2^*,\infty)$. 
\end{itemize}

A complete solution of case P is given in \cite{Bai10b}, where it was proved that the optimal
policies are simple lump sum
dividend strategies. This paper also covers the
more complex case R, but without a complete solution. To be more concrete, a few definitions are needed.
With $g$ a canonical solution, define for $a_1<a_2$,
$$I(a_1,a_2,c)=\int_{a_1}^{a_2}(k-cg'(x))dx= k(a_2-a_1) -c(g(a_2)-g(a_1)).$$
Define the set $C$, possibly empty, by $C=C_1\cup C_2$ where
\begin{eqnarray*}
C_1&=&\{c>0: \mbox{ there exists }0< \underline{u}<\bar{u} \mbox{ so that }
cg'(\underline{u})=cg'(\bar{u})=k \mbox{ and } I(\underline{u},\bar{u},c)=K\},\\
C_2&=&\{c>0: \mbox{ there exists } \bar{u}>0 \mbox{ so that }
cg'(\bar{u})=k,\;cg'(0)<k \mbox{ and } I(0,\bar{u},c)=K\}.
\end{eqnarray*}
Since $I$ and $g'$ are continuous, $C$ is closed. Hence if $C\neq \emptyset$, 
$c^*=\max\{c:c\in C\}$ is well defined. However, for known $c^*$, the corresponding $\underline{u}$
and $\bar{u}$ may not be unique. Therefore we define
$$U=\{u:c^*g'(u)=k\}.$$
Let
$$\bar{u}^*=\max\{u\in U: I(\underline{u},u,c^*)=K \mbox{ for some } \underline{u}\in \{0\}\cup U\}$$
and then
$$\underline{u}^*=\max\{u\in  \{0\}\cup U: I(u,\bar{u}^*,c^*)=K\}.$$

The case R can be divided into several subclasses as follows: 
\begin{itemize}
\item[R1.] $0<x_1^*<\underline{u}^*<x_2^*<\bar{u}^*,\;\;c^*g'(0)\geq k$.
\item[R2.] $0<x_1^*<\underline{u}^*<x_2^*<\bar{u}^*,\;\;c^*g'(0)<k$.
\item[R3.] $0=\underline{u}^*<x_1^*<x_2^*<\bar{u}^*,\;\;c^*g'(x_1^*)\geq k$.
\item[R4.] $0=\underline{u}^*<x_1^*<x_2^*<\bar{u}^*,\;\;c^*g'(x_1^*)< k$.
\item[R5.] $0=\underline{u}^*<\bar{u}^*<x_1^*,\;\;x_2^*<\infty,\;\;c^*g'(x_2^*)\geq k$.
\item[R6.] $0=\underline{u}^*<\bar{u}^*<x_1^*,\;\;x_2^*<\infty,\;\;c^*g'(x_2^*)< k$.
\item[R7.] $x_2^*=\infty$, i.e. $g$ is concave on $[x_1^*,\infty)$.
\item[R8.] None of the above.
\end{itemize}
As pointed out in Remark 2.2 in \cite{Bai10b}, the case R8 is pretty odd, so we drop it in
this paper.

If it exists, let $(c^*,\bar{u}^*)$ be a pair that satisfies
\begin{equation}
c^*g'(\bar{u}^*)=k\;\;\;\mbox{and}\;\;\;c^*g(\bar{u}^*)=k\bar{u}^*-K,
\label{22}\end{equation}
so that in particular
$$c^*=\frac{k}{g'(\bar{u}^*)}.$$ 
Under R5 and R6, $(c^*,\bar{u}^*)$ always exists, and is as given in the definition of R5-R6 above.   
 Under R7, if it does not exist, we can set $c^*=\bar{u}^*=0$. 

By the properties of $g$, $g'$ is ultimately increasing or decreasing, so we can define
$$g_{\infty}'=\lim_{x\rightarrow\infty}g'(x).$$
Under R1-R6, typically $g_{\infty}'=\infty$, while under R7, $g_{\infty}'$
is always finite. Also define
$$c^{\infty}=\frac{k}{g_{\infty}'}\;\;\mbox{with}\;\;c^{\infty}=0\;\;\mbox{if}\;\;
g_{\infty}'=\infty.$$

It was proved in \cite{Bai10b} that for cases R1-R4 the optimal policies are always
the simple barrier policies $\pi_{\bar{u}^*,\underline{u}^*}$ with corresponding value functions
$V_{\bar{u}^*,\underline{u}^*}(x)$.
It was  proved that this is also the case if
\begin{equation}
\mu(x_\lambda)k-\lambda(kx_\lambda -K)\leq 0
\label{11}\end{equation}
and R5 or R6 apply, and also if R7 applies,  (\ref{11}) holds and
 $c^*\geq c^\infty$. When R7 applies and
$c^*< c^\infty$ it was proved that there is no optimal policy, but
the value function equals $V^*(x)=c^{\infty} g(x)$.

Therefore, what remains is R5-R7 when

\begin{equation}
\mu(x_\lambda)k-\lambda(kx_\lambda -K)>0.
\label{23}\end{equation}
 
Figure 1.1 gives a graphical illustration of
cases R5 and R6. Note that in case R6, $K_3\leq \min\{K_2,K\}$, since otherwise it is possible to increase
$c^*$, bringing us to one of the cases R1-R4. Also if $K_3=K_2=K$, by definition of $c^*$, we have 
case R3. Therefore $K_3<\min\{K_2,K\}$.  But then $c^*\geq c^{\infty}$, or equivalently $g_{\infty}'\geq g'(\bar{u}^*)$, for otherwise
$K_3=\infty$. 

\begin{figure}[t]
\epsfig{file=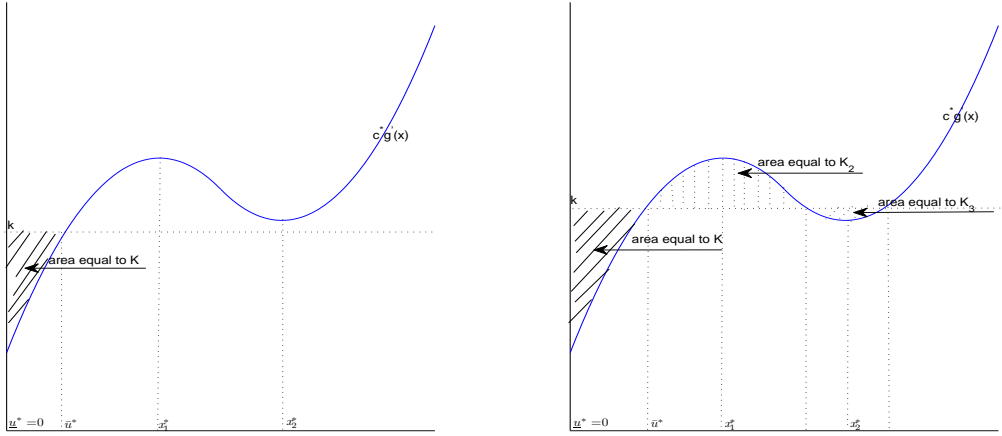,height=7cm,width=17cm}
\caption{Plots of cases R5 (left) and R6 (right). In case R6, $K_3<\min\{K_2,K\}$.}
\label{fig3}
\end{figure}

We can now formulate the two problems of this paper. A solution of these problems together with the results
in \cite{Bai10b} will give a complete solution to the dividend problem under A1-A4, with the exception of the
odd case R8.

\vspace{0.2cm}
\b{Problem 1.} Assume A1-A4 and that R5 or R6 apply so that in particular $c^*\geq c^\infty$. Also assume
(\ref{23}). Find the optimal value function, and if it exists, the optimal dividend strategy.

\vspace{0.2cm}
\b{Problem 2.} Assume A1-A4 and that R7 applies. Also assume  (\ref{23})  and that 
$c^*\geq c^\infty$. Find the optimal value function, and if it exists, the optimal dividend strategy.
\vspace{0.2cm}

It was proved in \cite{Bai10b} that a simple
lump sum dividend strategy cannot be optimal  for any of these two problems. 
This kind of result is not new for diffusion processes, see for example
Example 4.3 in \cite{Eg08}, and for the case with no fixed transaction costs, see \cite{Al98}.
For more background information and details, the reader should consult \cite{Bai10b} and references therein.

As pointed out in Remark 2.3 in \cite{Bai10b}, in order for R5 or R6 to apply it is
necessary that $\mu(0)<-\frac{\lambda K}{k}$. It can thus be argued that these cases
are less interesting from a practical point of view. However, from a theoretical
point of view these cases are the most interesting, and certainly the most challenging,
and as we shall see in Section 3, the optimal solution is  highly nontrivial. 

A complete solution of Problem 1 is given in Section 3, and of Problem 2 in Section 4. However, before we can 
present the solutions, several definitions and preliminary results are needed. This is the topic of the next section.

\newsection{Some preliminary results}
 The notation and definitions are the same as in
Section 1. Throughout this section, the assumptions stated in Problem 1 or Problem 2 are
assumed to hold, so in particular (\ref{22}) is assumed to have a solution with $c^*>0$.
All results  are proved in the appendix.

The following lemma yields two solutions  of $Lg(x)=0$ that are useful for 
further analysis. 

\begin{lemma} There exists two independent, positive
solutions $g_{1}$ and $g_{2}$ of $Lg(x)=0$ such that:
\begin{itemize}
\item[(i)] $g_{1}(0)=1,$
$g_{1}'(x)<0$ on $[0,x_{\lambda})$, $g_{1}'(x_{\lambda})=0,$ $g_{1}'(x)>0$ on $(x_{\lambda}, \infty)$, and $g_{1}''(x)>0$
on $[x_{\lambda}, \infty).$
Furthermore, there exists at most one point $z_0>0$  such that $g_{1}''(z_{0})=0$,
and at this point $g'''(z_0)>0$. If $z_0$ exists then $z_0\in (0,x_1^*)$. Also $g_1(x)>0$ for all $x\geq 0$.
\item[(ii)] $g_{2}$ is strongly increasing on $[0,\infty)$ with $g_2(0)=1$ and $g_2'(0)>0$.
\end{itemize}
\end{lemma}

It follows from Lemma 2.1 that such a $z_0$ exists if and only if
$g_1''(0)<0$. 
If the function $g_1$ in Lemma 2.1(i) satisfies  $g_1''(0)\geq 0$ we set $z_0=0$.
In this case, the interval $(0,z_0)$ is just the emptyset.

With $g_1$ and $g_2$ as in Lemma 2.1, a canonical solution becomes
$$g(x)=g_2(x)-g_1(x).$$
In the sequel, the canonical solution $g(x)$ will always be this particular function.
Define
$$g(x;\beta)=g_2(x)-\beta g_1(x),$$
so that $g(x;1)=g(x)$, the canonical solution.
When writing $g'(x;\beta)$ we shall mean
$\frac{d}{dx}g(x;\beta)=g_2'(x)-\beta g_1'(x)$. Similarly
with $g''(x;\beta)$ and $g'''(x;\beta)$. 

From the fact that $g'(x)>0$ and $g_2'(x)>0$ it follows easily that $g'(x;\beta)>0$
for  $\beta\in [0,1]$, hence we can define
$$\gamma(\beta,x)=k\left(\frac{g(x;\beta)}{g'(x;\beta)}-x\right)+K.$$
Note that 
\begin{equation}
\gamma(1,\bar{u}^*)=0.
\label{21}\end{equation}

The following function will play an important role in this paper.
\begin{equation}
h(x)=\frac{g_2''(x)}{g_1''(x)},\;\;\;\;x\neq z_0.\label{20}
\end{equation}
Here is a simple observation that will often be used
\begin{eqnarray}
g''(x;\beta)=g_1''(x)(h(x)-\beta),\;\;\;\;x\neq z_0.\label{26}
\end{eqnarray}

The next result contains some important properties of the function $h$.
\begin{lemma}
The function $h$ is  strongly decreasing on $(0,z_0)\cup (z_0,x_\lambda)$ and strongly
increasing on $(x_\lambda,\infty)$. Thus the limit $h_\infty=\lim_{x\rightarrow\infty}h(x)$ exists, but may be infinite. 
In cases R5 and R6, $h_\infty\geq 1$, while in case R7, $h_\infty \leq 1$.
\end{lemma}

Using that $z_0<x_1^*$ and the fact that $h(x_1^*)=1$ (see Lemma A.1(c)),
it follows from Lemma 2.2 that we can define a continuously differentiable, strongly decreasing function $h_-^{-1}:[h(x_\lambda),1]$
onto  $[x_1^*,x_\lambda]$ so that $h(h_-^{-1}(x))=x$. Furthermore, we
can define a continuously differentiable, strongly increasing function $h_+^{-1}:[h(x_\lambda),h_\infty)\rightarrow [x_\lambda,\infty)$
so that $h(h_+^{-1}(x))=x$. See Figure 2.1.

\begin{figure}[t]
\epsfig{file=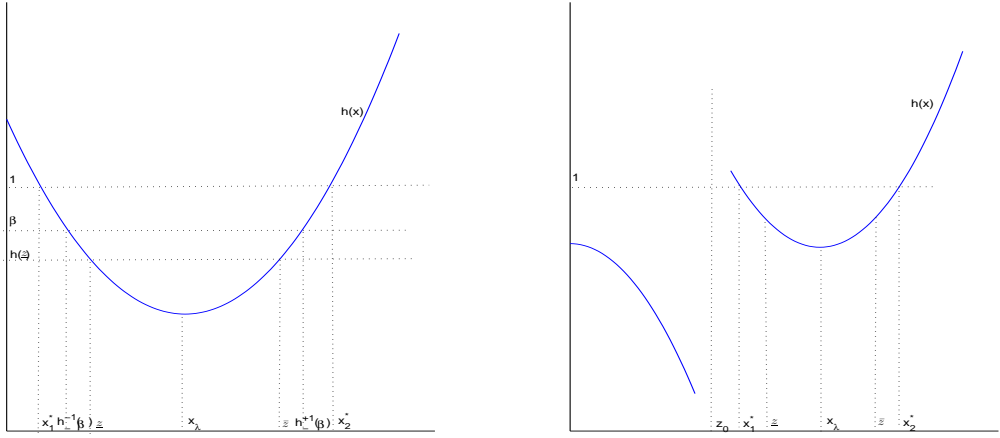,height=7cm,width=17cm}
\caption{Plots of $h(x)$. In the left plot, $g_1''(0)\geq 0$, while
in the right plot $g_1''(0)<0$.}
\label{fig1}
\end{figure}

By Lemma 2.1 in \cite{Bai10b} there is a unique $\underline{z}\in (x_1^*,x_\lambda)$
so that
\begin{equation}\mu(\underline{z})k-\lambda(k\underline{z}-K)=0\label{24}\end{equation}
and
\begin{equation}\mu(x)k-\lambda(kx-K)<0,\;\;\;\;x<\underline{z}.
\label{25}\end{equation}
In Lemma 2.1 in \cite{Bai10b} it is assumed that R5 or R6 apply, but looking at the
proof shows that the only requirement is that $\bar{u}^*<x_1^*$. Therefore, (\ref{24}) and
(\ref{25}) are valid under R7 as well whenever (\ref{22}) holds.

\begin{lemma} Let $\underline{z}$ be as in (\ref{24}). Then
$h(\underline{z})<1$. Furthermore,
if $g_1''(0)<0$, or equivalently if  $z_0>0$, then 
$h(0)\leq h(\underline{z})$. Finally, for $\beta\in (h(\underline{z}),1]$, 
$$g''(x;\beta)>0,\;\;\;\;x\in (0,h_-^{-1}(\beta)).$$
\end{lemma}
By Lemmas 2.2 and 2.3,  $h(\underline{z})<h_\infty$ in cases R5 and R6, while
in case R7 this inequality may not hold.

\begin{lemma} 
There is a unique, continuously differentiable, strongly decreasing
function $u_1: [h(\underline{z}),1]\rightarrow [\bar{u}^*,\underline{z}]$   so that
$$\gamma(\beta,u_1(\beta))=0.$$
Furthermore, $\bar{u}^*<u_1(\beta)<h_-^{-1}(\beta)$ for $\beta\in (h(\underline{z}),1)$, while
$u_1(h(\underline{z}))=\underline{z}$ and $u_1(1)=\bar{u}^*$.
\end{lemma}  

\vspace{0.2cm}

It follows from Lemma 2.3 that for $\beta\in [h(\underline{z}),1]$, as a function
of $x$, $g'(x;\beta)$ increases on $(0,h_-^{-1}(\beta))$. Furthermore, if $h(\underline{z})<h_\infty$ and
$\beta\in [h(\underline{z}),h_\infty)$, it follows from  (\ref{26}) 
and the fact that $z_0<x_1^*<h_-^{-1}(\beta)$,
that $g'(x;\beta)$  
decreases on $(h_-^{-1}(\beta),h_+^{-1}(\beta))$ and
increases on $(h_+^{-1}(\beta),\infty)$. Therefore, when $h(\underline{z})<h_\infty$ we can define 
\begin{equation}
v(\beta)=\left\{\begin{array}{lll}\frac{g'(h_+^{-1}(\beta);\beta)}{g'(u_1(\beta);\beta)},
\;\;\;&\beta\in [h(\underline{z}),1]\cap [h(\underline{z}),h_\infty),\\
\infty,&\mbox{otherwise}. \end{array}\right.
\label{27}
\end{equation}
Now $u_1(h(\underline{z}))=\underline{z}=h_-^{-1}(h(\underline{z}))$, so  we get
since $g'(x;\beta)$ has a local maximum at $h_-^{-1}(\beta)$,
$$v(h(\underline{z}))=\frac{g'(h_+^{-1}(h(\underline{z}));h(\underline{z}))}{g'(h_-^{-1}(h(\underline{z}));h(\underline{z}))}<1.$$
Therefore, if  $h(\underline{z})<h_\infty$, define
\begin{equation}
\beta_0=\sup\{\beta_1\in (h(\underline{z}),1]: v(\beta)\leq 1\;\;\mbox{for all}\;\;\beta\in (h(\underline{z}),\beta_1]\}.
\label{28}\end{equation}
In case R7, $h_\infty\leq 1$ by Lemma 2.2, so then $h_+^{-1}(1)$ is not defined, yielding that $\beta_0<1$.
In cases R5 and R6,
$$v(1)=\frac{g'(h_+^{-1}(1))}{g'(u_1(1))}=\frac{g'(x_2^*)}{g'(\bar{u}^*)}=\frac{c^*}{k} g'(x_2^*).$$
Therefore, when there is a strict inequality in R5, $\beta_0<1$.

To prepare for Problem 2, assume that $g(x)$ is concave for $x>x_\lambda$. Also,
for the moment we do not assume that $h(\underline{z})<h_\infty$. 
By Lemma 2.2,  $h_{\infty}\leq 1$. 
Since $g$ is ultimately concave, $g_\infty'<\infty$
and therefore
$$g_{1,\infty}':=\lim_{x\rightarrow\infty}g_1'(x)<\infty\;\;\;\mbox{if and only if}
\;\;\;g_{2,\infty}':=\lim_{x\rightarrow\infty}g_2'(x)<\infty.$$
Assume that $g_{1,\infty}'<\infty$. By Lemma 2.4 we can define
$$G(\beta)=\frac{ g_{2,\infty}'-\beta g_{1,\infty}'}{g'(u_1(\beta);\beta)}
=\frac{\lim_{x\rightarrow\infty}g'(x;\beta)}{g'(u_1(\beta);\beta)},\;\;\;\;\beta\in [h(\underline{z}),1].$$
Note that  by Lemmas 2.1, 2.2 and 2.4, $g'(u_1(\beta);\beta)>g'(0,\beta)=1-\beta\geq 0$ for $\beta\in [h(\underline{z}),1].$

\begin{lemma} Assume that $g$ is concave
for $x>x_\lambda$ and that  $g_{1,\infty}'<\infty$. Then the function $G$ is
continuously differentiable and strongly increasing on $(h(\underline{z}),1)$.
\end{lemma}

Let us return to the general case and again assume that $h(\underline{z})<h_\infty$. Since $g'(x;\beta)$ increases
for $x> h_+^{-1}(\beta)$,
$\lim_{x\rightarrow\infty}g'(x;\beta)$ always exists, but may
be infinite. In that case we set $G(\beta)=\infty$.
Again by the ultimate increase of $g'(x;\beta)$,  $\lim_{x\rightarrow\infty}g'(x;\beta)>g'(h_+^{-1}(\beta);\beta),$ hence
\begin{equation}
G(\beta)>v(\beta),\;\;\;\;\beta \in [h(\underline{z}),\beta_0).\label{29}
\end{equation}
Finally, if it exists, define $\alpha_0$ by
\begin{equation}
\alpha_0=\left\{\begin{array}{lll}h(\underline{z})\;\;\;&\mbox{if}\;\;G(h(\underline{z}))\geq 1,\\
\alpha_1  \;\;\;&\mbox{if}\;\;G(h(\underline{z}))<1,\end{array}\right.
\label{210}\end{equation}
where $\alpha_1\in (h(\underline{z}),\beta_0)$ is the unique value 
that satisfies $G(\alpha_1)=1$. So by definition, if it exists, $\alpha_0<\beta_0$.
By (\ref{29}), such an $\alpha_0$ exists whenever $\beta_0<1$. When $\beta_0=1$,
$\alpha_0$ exists if and only if $G(1)>1$, i.e. if and only if $g_{\infty}'>g'(\bar{u}^*)$,
or equivalently if and only if $c^*>c^{\infty}$. In this case, by definition of $c^*$, (\ref{22}) will
necessarily hold.

\vspace{0.2cm}
\begin{lemma}
Assume that $h(\underline{z})<h_\infty$ and that $\alpha_0$ in (\ref{210}) exists.
Then there are continuously differentiable functions $u_i$,
i=1,2,3, defined on $(\alpha_0,\beta_0)$, so that $\gamma(\beta,u_1(\beta))=0$ and
\begin{equation}
g'(u_1(\beta);\beta)=g'(u_2(\beta);\beta)= g'(u_3(\beta);\beta).\label{211}
\end{equation}
Furthermore, $u_1$ is strongly decreasing and
$$\bar{u}^*<u_1(\beta)<h_-^{-1}(\beta)<\underline{z}<u_2(\beta)<h_+^{-1}(\beta)<u_3(\beta).$$
Also if $\alpha_0=h(\underline{z})$,
\begin{equation}
\lim_{\beta\downarrow \alpha_0}u_1(\beta)=\lim_{\beta\downarrow \alpha_0}u_2(\beta)=\underline{z},
\label{212}\end{equation}
and if $\alpha_0>h(\underline{z})$,
\begin{equation}
\lim_{\beta\downarrow \alpha_0}u_3(\beta)=\infty. \label{213}
\end{equation}
Finally, if $\beta_0<\min\{1,h_\infty\}$.
\begin{equation}
\lim_{\beta\uparrow \beta_0}u_2(\beta)=\lim_{\beta\uparrow \beta_0}u_3(\beta).
\label{214}\end{equation}
\end{lemma}

Define the function $J$  by
$$J(\beta,u_1,u_2)=k\int_{u_1}^{u_2}\left(1-\frac{g'(x;\beta)}{g'(u_2;\beta)}\right)dx
=k\left(u_2-u_1-\frac{g(u_2;\beta)-g(u_1;\beta)}{g'(u_2;\beta)}\right).$$
When (\ref{211}) is satisfied, we will use the simplified notation
\begin{eqnarray*}
J_1(\beta)&=&J(\beta,u_1(\beta),u_2(\beta)),\\
J_2(\beta)&=&J(\beta,u_2(\beta),u_3(\beta)),\\
J_{13}(\beta)&=&J(\beta,u_1(\beta),u_3(\beta)),
\end{eqnarray*}
so that in particular $J_{13}(\beta)=J_1(\beta)+J_2(\beta)$.
\begin{lemma} Under the assumptions of Lemma 2.6, for $\beta \in (\alpha_0,\beta_0)$, both
$J_1$ and $J_{13}$ are strongly decreasing.
\end{lemma}

When $G(h(\underline{z}))\geq 1$ then $\alpha_0=h(\underline{z})$ and so
$J_1(\alpha_0)=\lim_{\beta\downarrow \alpha_0}J_1(\beta)=0$ by (\ref{212}). When
$G(h(\underline{z}))<1$, $\alpha_0>h(\underline{z})$ and then
$J_1(\alpha_0)<0$. By Lemma 2.6, 
$$J_2(\alpha_0)=\lim_{\beta\downarrow \alpha_0}J_{13}(\beta)
-\lim_{\beta\downarrow \alpha_0}J_1(\beta)$$
exists, but may be infinite because of (\ref{213}). To get the precise solutions of our problems 
it is convenient to split into five different cases. Again, the assumption $h(\underline{z})<h_\infty$ is
in force.

\begin{itemize}
\item[G1.] $G(h(\underline{z}))\geq 1$.
\item[G2.] $G(h(\underline{z}))< 1$, $\alpha_0$ exists  and $J_{13}(\alpha_0)\geq 0$.
\item[G3.] $G(h(\underline{z}))< 1$, $\alpha_0$ exists,  $J_{13}(\alpha_0)<0$ and $J_2(\alpha_0)>K$.
\item[G4.] $G(h(\underline{z}))< 1$, $\alpha_0$ exists, $J_{13}(\alpha_0)<0$ and $J_2(\alpha_0)\leq K$.
\item[G5.] $\alpha_0$ does not exist, that is $c^*\leq c^{\infty}$.
\end{itemize}
In G1, $\alpha_0$ exists by definition. Also, in Problem 1, G5 becomes $c^*=c^{\infty}$.

\vspace{0.2cm}
\b{Remark.} Under R5 or R6
the assumption that $G(h(\underline{z}))>1$, or equivalently that
$\lim_{x\rightarrow\infty}g'(x;h(\underline{z}))>g'(\underline{z};h(\underline{z}))$
is extremely weak. 
For these cases, typically $g_{\infty}'=\infty$, and then
$$\lim_{x\rightarrow\infty}g'(x;h(\underline{z}))=g_{\infty}'
+(1-h(\underline{z}))\lim_{x\rightarrow\infty}g_1'(x)=\infty.$$
It is proved in \cite{Bai11}, Proposition 2.5, that a sufficient condition for $g_{\infty}'=\infty$ is that there exist
an $x_0>0$ and an $\varepsilon>0$ so that
$\mu'(x)<\lambda-\varepsilon$ for all $x>x_0$.

\vs
The next result is the basis for the optimality result Theorem 3.1 for Problem 1.

\begin{proposition}
Given the assumptions of Problem 1  and either of G1, G2 or G3.
Then 
there exists a  $\tilde{\beta}\in(\alpha_0,\beta_0)$ and numbers
$\tilde{u}_i=u_i(\tilde{\beta}),\;i=1,2,3,$ with
$$\bar{u}^*< \tilde{u}_1<h_-^{-1}(\tilde{\beta})<\underline{z}<\tilde{u}_2
<h_+^{-1}(\tilde{\beta})<\tilde{u}_3,$$
so that
\begin{itemize}
\item[(i)] $\gamma(\tilde{\beta},\tilde{u}_1)=0$.
\item[(ii)] $g'(\tilde{u}_1;\tilde{\beta})=g'(\tilde{u}_2;\tilde{\beta})
=g'(\tilde{u}_3;\tilde{\beta})$.
\item[(iii)] exactly one of the following two possibilities holds:
\begin{itemize}
\item[(iii-a)] $-J(\tilde{\beta},\tilde{u}_1,\tilde{u}_2)=
J(\tilde{\beta},\tilde{u}_2,\tilde{u}_3)\leq K$.
\item[(iii-b)] $-J(\tilde{\beta},\tilde{u}_1,\tilde{u}_2)>
J(\tilde{\beta},\tilde{u}_2,\tilde{u}_3)= K$.
\end{itemize}
\end{itemize}
In case of G3, (iii-b) applies.
\end{proposition}

Similarly to Proposition 2.1, the next result is the basis for Theorem 4.1 that partially covers
Problem 2.
\begin{proposition}
Given the assumptions of Problem 2 and  either of G1, G2 or G3.
Also assume  one of the following two conditions:
\begin{enumerate}
\item $h_{\infty}=1$.
\item $h(\underline{z})<h_{\infty}<1$ and $G(h_\infty)>1$. 
\end{enumerate}
Then
the results in Proposition 2.1 still hold.
\end{proposition}

\newsection{Solution of Problem 1}
 The notation and definitions are  the same as in
Sections 1 and 2. Throughout this section, the assumptions stated in Problem 1  are
assumed to hold.
All results  are proved in the appendix.

In the presentation of the results, we need one more definition. Let $\gamma\in [h(\underline{z}),1]$ and define
$$c_{\gamma}=\frac{k}{g'(u_1(\gamma);\gamma)},$$
so in particular $c_1=c^*$ if the latter is positive.

\vspace{0.2cm}

\begin{theorem} 
Given the same assumptions and notation as in Proposition 2.1:
\begin{itemize}
\item[(a)] If (i), (ii) and (iii-a) of Proposition 2.1 hold,
then the two-level lump sum dividend strategy 
$\pi_{(\bar{u}^*,0)(\tilde{u}_1,\tilde{u}_3,0)}$ is optimal. The corresponding 
value function is given as
$$V_{(\bar{u}^*,0)(\tilde{u}_1,\tilde{u}_3,0)}(x)=
\left\{\begin{array}{llll}
c^*g(x),\;\;\;&x\in[0,\bar{u}^*),\\
kx-K,&x\in [\bar{u}^*,\tilde{u}_1),\\
c_{\tilde{\beta}}g(x;\tilde{\beta}),& x\in [\tilde{u}_1,\tilde{u}_3),\\
kx-K,&x\in [\tilde{u}_3,\infty).
\end{array}\right.$$
\item[(b)] If (i), (ii) and (iii-b) of Proposition 2.1 hold,
then the two-level lump sum dividend strategy 
$\pi_{(\bar{u}^*,0)(\tilde{u}_1,\tilde{u}_3,\tilde{u}_2)}$ is optimal. The corresponding 
value function is given as
$$V_{(\bar{u}^*,0)(\tilde{u}_1,\tilde{u}_3,\tilde{u}_2)}(x)=
\left\{\begin{array}{llll}
c^*g(x),\;\;\;&x\in[0,\bar{u}^*),\\
kx-K,&x\in [\bar{u}^*,\tilde{u}_1),\\
c_{\tilde{\beta}}g(x;\tilde{\beta}) ,& x\in [\tilde{u}_1,\tilde{u}_3),\\
c_{\tilde{\beta}}g(\tilde{u}_2;\tilde{\beta})+k(x-\tilde{u}_2)-K,&x\in [\tilde{u}_3,\infty).
\end{array}\right.$$
\end{itemize}
In both cases, the value function is continuously differentiable.
\end{theorem}

\vspace{0.3cm}
When Theorem 3.1(a) holds, the first payment will lead to ruin, while when
Theorem 3.1(b) holds, the first payment results in ruin only if
$X_{\tau_1^{\pi^*}}^{\pi^*}\in [\bar{u}_1^*,\tilde{u}_1]$. 

\vspace{0.3cm}
Now turn to the case G4, and here the solution is  different from
that in Theorem 3.1.

\begin{theorem} 
Assume that case G4 applies.
Then the optimal value function is given as
\begin{eqnarray*}
V^*(x)=
\left\{\begin{array}{llll}
c^*g(x),\;\;\;&x\in[0,\bar{u}^*),\\
kx-K,&x\in [\bar{u}^*,u_1(\alpha_0)),\\
c_{\alpha_0}g(x;\alpha_0),& x\in [u_1(\alpha_0),\infty).\\
\end{array}\right.
\end{eqnarray*}
When $x\in (0,u_1(\alpha_0)]$ the lump sum dividend barrier strategy 
$\pi_{\bar{u}^*,0}$ is optimal, while for $x>u_1(\alpha_0)$
 there is no optimal strategy and $V^*(x)={\mathop{\lim}_{\bar{u}\rightarrow
\infty}V_{(\bar{u},0)}(x)}$.
\end{theorem}

Finally, we give the  solution for G5.

\begin{theorem} 
Assume that case G5 applies so that $c^*=c^{\infty}$. Then there is no
optimal policy, but the value function is given as
$$V^*(x)=c^*g(x),\;\;\;\;x\geq 0.$$
\end{theorem}

The result in Theorem 3.3 is just a limit of those in Theorems 3.1 and 3.2. For
example in Theorem 3.2, since as
$c^{\infty}\uparrow c^*$, $G(1)\downarrow 1$ and so $\alpha_0\uparrow 1$. But then
$u_1(\alpha_0)\downarrow u_1(1)=\bar{u}^*$ and by (\ref{a4}), $c_{\alpha_0}\uparrow c^*$.
Hence the interval $[\bar{u}^*,u_1(\alpha_0))\rightarrow \emptyset$ and 
$c_{\alpha_0}g(x;\alpha_0)\rightarrow c^*g(x)$.
Therefore
$$V^*(x)\rightarrow c^*g(x).$$
Similarly, in Theorem 3.1 as $c^{\infty}\uparrow c^*$, $\tilde{\beta}\uparrow 1$, $\tilde{u}_1
\downarrow \bar{u}^*$, $c_{\tilde{\beta}}\uparrow c_1=c^*$ and $\tilde{u}_3\rightarrow \infty$.

\newsection{Solution of Problem 2}
 The notation and definitions are  again the same as in
Sections 1 and 2. Throughout the section, the assumptions stated in Problem 2  are
assumed to hold.
All results  are proved in the appendix.

As in the last section, we start with the case $c^*>c^{\infty}$, or equivalently 
$g_{\infty}'> g'(\bar{u}^*)$. The case $c^*=c^{\infty}$
is then solved in Theorem 4.4.

\begin{theorem}
Under the assumptions and notation of Proposition 2.2, the optimal policy and value function
are as in Theorem 3.1.
\end{theorem}

It remains to solve the cases not covered by Theorem 4.1, i.e. when either
\begin{itemize}
\item[H1.] $h(\underline{z})\geq h_\infty$.
\item[H2.] $h(\underline{z})<h_\infty<1$ and $G(h_\infty)\leq 1$.
\item[H3.] Case G4.
\item[H4.] $c^*=c^{\infty}$.
\end{itemize}
Note that there can be  overlap between case H3 and either of H1 or H2.
For H1, since $h(\underline{z})<h(x_1^*)=1$, the condition $h_\infty <1$
is automatically satisfied.

\begin{lemma} Assume that  one of H1-H3 apply.
Then there is a unique $\hat{\beta} \in 
(h(\underline{z}),1)$ that satisfies $G(\hat{\beta})=1$. Furthermore, $u_1(\hat{\beta})<\underline{z}$.
\end{lemma}

\begin{theorem} 
Assume that  one of H1-H3 apply, and let $\hat{\beta}$ be as in Lemma 4.2.
Then the optimal value function is given as
\begin{eqnarray*}
V^*(x)=
\left\{\begin{array}{llll}
c^*g(x),\;\;\;&x\in[0,\bar{u}^*),\\
kx-K,&x\in [\bar{u}^*,u_1(\hat{\beta})),\\
c_{\hat{\beta}}g(x;\hat{\beta}),& x\in [u_1(\hat{\beta}),\infty).\\
\end{array}\right.
\end{eqnarray*}
When $x\in (0,u_1(\hat{\beta})]$ the lump sum dividend barrier strategy 
$\pi_{\bar{u}^*,0}$ is optimal, while for $x>u_1(\hat{\beta})$
 there is no optimal strategy and $V^*(x)={\mathop{\lim}_{\bar{u}\rightarrow
\infty}V_{(\bar{u},0)}(x)}$.
\end{theorem}

Finally, here is the result when $c^*=c^{\infty}$.

\begin{theorem} 
Assume H4, i.e. that $c^*=c^{\infty}$. Then there is no
optimal policy, but the value function is given as
$$V^*(x)=c^*g(x),\;\;\;\;x\geq 0.$$
\end{theorem}
Just as pointed out after Theorem 3.3, Theorem 4.4 is a limiting case of Theorem 4.3 as
$c^{\infty}\uparrow c^*$.

\appendix
\newsection{Appendix: Proofs }
\b{Proof of Lemma 2.1.} The proof that $z_0<x_1^*$ is deferred
until after the proof of Lemma A.2. For the construction of the $g_i$,
it is well known, see e.g. \cite{Bor},
that the equation $Lg(x)=0$ has a strongly increasing as well as a strongly decreasing 
solution, both positive. Let $g_2$ be the strongly increasing solution, divided 
by $g_2(0)$ so that the new $g_2$  satisfies $g_2(0)=1$.

Using the same arguments as in the proof of Lemma A.5 in \cite{P4}, but starting 
at $x_\lambda$, we find that there is a solution $g_1$ satisfying
$$g_1(x_\lambda)=1,\;\;g_1'(x_\lambda)=0,\;\;g_1'(x)>0\;\mbox{on}\;
(x_\lambda,\infty)\;\;\mbox{and}\;\;g_1''(x)>0\;\mbox{on}\;
[x_\lambda,\infty).$$
We will show that $g_1'(x)<0$ on $(0,x_\lambda)$, which clearly implies that $g_1$
is positive.
Since $g_1'(x_\lambda)=0$ and $g_1''(x_\lambda)>0$, there is an
$\varepsilon>0$ so that $g_1'(x)<0$ and $g_1''(x)>0$ on $(x_\lambda-\varepsilon,
x_\lambda)$. Define, if it exists
$$x_0=\max\{x\leq x_\lambda-\varepsilon: g_1'(x)=0\}.$$
Since $g_1'(x)<0$ for $x\in (x_0,x_\lambda)$, this implies that
$g_1''(x_0)\leq 0$ and that $g_1(x_0)>1$. The relation $Lg_1(x)=0$ gives
\begin{equation}
g_1''(x)=\frac{2\lambda}{\sigma^2(x)}g_1(x)-\frac{2\mu(x)}{
\sigma^2(x)}g_1'(x). \label{a1}\end{equation}
Therefore, by the definition of $x_0$, 
$g_1''(x_0)= \frac{2\lambda}{\sigma^2(x)}g_1(x_0)>0$, leading to
a contradiction. Consequently, $g_1'(x)<0$ on $(0,x_\lambda)$. 

We now turn to the proof that there is at most one $z_0<x_\lambda$ so that $g_1''(z_0)=0$.
Differentiating $Lg_1(x)=0$ gives
\begin{equation}
g_1'''(x)=\frac{2(\lambda-\mu'(x))}{\sigma^2(x)}g_1'(x)-\frac{2}{\sigma^2(x)}
(\mu(x)+\sigma(x)\sigma'(x))g_1''(x). \label{a2}\end{equation}
If there is a $z_0<x_\lambda$ so that $g_1''(z_0)=0$, it follows from
(\ref{a2}) that $g_1'''(z_0)>0$. But then clearly there can only be one such $z_0$.
Also $\min_{x\geq 0}g_1(x)=g_1(x_\lambda)=1$ so $g_1(x)>0$ for all $x\geq 0$.
Finally, we rescale $g_1$ by dividing it by $g_1(0)$ so that the new $g_1$ satisfies
$g_1(0)=1$.

\vspace{0.2cm}

Here are some simple results that will be useful in the sequel.

\begin{lemma}
Let $g_1$ and $g_2$ be as in Lemma 2.1, and let $h$ be as in (\ref{20}). 
\begin{itemize}
\item[(a)] The Wronskian $W(x)=g_1(x)g_2'(x)-g_2(x)g_1'(x)>0$.
\item[(b)] $\frac{g_2(x)}{g_1(x)}$ is increasing in $x$.
\item[(c)] $h(x_1^*)=1$ and if $x_2^*<\infty$, $h(x_2^*)=1.$
\item[(d)] $$h'(x)=\frac{4\lambda(\lambda-\mu'(x))W(x)}
{\sigma^4(x)(g_1''(x))^2},\;\;\;x\neq z_0.$$
\item[(e)] $$\frac{g(x;h(x))}{g'(x;h(x))}=\frac{\mu(x)}{\lambda},\;\;\;x\neq z_0.$$
\item[(f)] $$\frac{\partial}{\partial \beta}\gamma(\beta,x)
=-k\frac{W(x)}{(g'(x;\beta))^2}<0.$$
\item[(g)] $$\frac{\partial}{\partial x}\gamma(\beta,x)
=-k \frac{g(x;\beta)}{(g'(x;\beta))^2}g_1''(x)(h(x)-\beta)
=-k \frac{g(x;\beta)}{(g'(x;\beta))^2}(g_2''(x)-\beta g_1''(x)).$$
\end{itemize}
\end{lemma}

{\sl Proof.} For (a),  we have that $W(0)=g'(0)>0$, and since
the Wronskian of two independent solutions never vanishes, the result follows.
As for (b), direct differentiation gives
$$\frac{d}{dx}\frac{g_2(x)}{g_1(x)}=\frac{W(x)}{(g_1(x))^2}>0$$
by (a).
For (c), we have
$$0=g''(x_1^*)=g_2''(x_1^*)-g_1''(x_1^*)\:\Rightarrow \;
h(x_1^*)=1.$$
Similarly, if $x_2^*<\infty$, $h(x_2^*)=1$.
For (d), differentiation gives
$$h'(x)=\frac{g_1''(x)g_2'''(x)-g_1'''(x)g_2''(x)}{(g_1''(x))^2}.$$
Using (\ref{a2}) for $g_i'''(x)$, $i=1,2$, yields
$$g_1''(x)g_2'''(x)-g_1'''(x)g_2''(x)=
\frac{2(\lambda-\mu'(x))}{\sigma^2(x)}(g_1''(x)g_2'(x)
-g_2''(x)g_1'(x)).$$
Using (\ref{a1}) for $g_i''(x)$, $i=1,2$, yields
\begin{equation}
g_1''(x)g_2'(x)-g_2''(x)g_1'(x)=\frac{2\lambda}{\sigma^2(x)}W(x).\label{a3}
\end{equation}
Combining these results proves (d). For (e), multiplying by $\frac 12 \sigma^2(x)
g_1''(x)$ in both the numerator and the denominator gives
\begin{eqnarray*}
\frac{g(x;h(x))}{g'(x;h(x))}&=&
\frac{\frac 12 \sigma^2(x)g_1''(x)g_2(x)-\frac 12 \sigma^2(x)g_2''(x)g_1(x)}
{\frac 12 \sigma^2(x)g_1''(x)g_2'(x)-\frac 12 \sigma^2(x)g_2''(x)g_1'(x)}\\
&=&\frac{(\lambda g_1(x)-\mu(x)g_1'(x))g_2(x)-(\lambda g_2(x)-\mu(x)g_2'(x))g_1(x)}
{(\lambda g_1(x)-\mu(x)g_1'(x))g_2'(x)-(\lambda g_2(x)-\mu(x)g_2'(x))g_1'(x)}\\
&=&\frac{\mu(x)W(x)}{\lambda W(x)}\\
&=&\frac{\mu(x)}{\lambda}.
\end{eqnarray*}
Finally, (f) and (g) follow by straightforward differentiation and use of the
definition of $h(x)$.

\vspace{0.2cm}
\b{Proof of Lemma 2.2.}
The first part is a direct consequence of Lemma 2.1(i) and Lemma A.1(a) and (d).
As for the second part, assume R5 or R6. Then for $x>x_2^*$, 
$0<g''(x)=g_2''(x)-g_1''(x)$. But since $g_1''(x)>0$ for $x>x_\lambda$,
it follows that $h(x)>1$ for $x>x_2^*$. This covers cases R5 and R6, and
the result for case R7 is proved in the same way.

\begin{lemma}
Assume (\ref{23}). Then 
\begin{itemize}
\item[(a)] $\gamma(h(x),x)<0,\;\;\;x<\underline{z}.$
\item[(b)] $\gamma(h(\underline{z}),\underline{z})=0.$
\end{itemize}
\end{lemma}

{\sl Proof.}
To prove (a), we use Lemma A.1(e)  together with (\ref{22}), so that for
$x<\underline{z}$,
$$\gamma(h(x),x)=k\frac{g(x;h(x))}{g'(x;h(x))}-kx+K
=\frac{1}{\lambda} (k\mu(x)-\lambda(kx-K))<0.$$
Part (b)  follows from 
the proof of (a) and the definition of $\underline{z}$.

\vspace{0.2cm}

\b{Proof that $z_0<x_1^*$ in Lemma 2.1.} By the comment right after Lemma 2.1,
$g_1''(0)<0$ and
since $g$ is convex on $[0,x_1^*]$, $0\leq g''(0)=
g_2''(0)-g_1''(0)$, i.e. $h(0)\leq 1$. By Lemma 2.1(i), $z_0<x_\lambda$,
so by Lemma A.1(d), $h$ is decreasing on $(0,z_0)$, in particular
$h(x)<1$ on $(0,z_0)$. But by definition of $z_0$, $g_1''(x)<0$ on $(0,z_0)$,
and so  $g''(x)=g_1''(x)(h(x)-1)>0$
on $(0,z_0)$. Consequently by continuity, $z_0\leq x_1^*$.
If $z_0=x_1^*$, then $h(x_1^*)=1$ implies that
$g_1''(z_0)=g_2''(z_0)=0$, which in turn implies that
$g_1''(z_0)g_2'(z_0)=
g_2''(z_0)g_1'(z_0).$
But by (\ref{a3}), this is equivalent to 
the Wronskian $W(z_0)=0$, a contradiction. Hence $z_0<x_1^*$.
Note that the proof is independent of (\ref{23}) and (\ref{25}).

\vspace{0.2cm}

\b{Proof of Lemma 2.3.} Since $\underline{z}\in (x_1^*,x_\lambda)$ and
$h$ decreases in this interval, the fact that $h(x_1^*)=1$ yields that $h(\underline{z})<1$.

Next assume that $z_0>0$, and also assume that $h(0)>h(\underline{z})$. We will show
that this gives a contradiction. 
Let $\beta\in
(h(\underline{z}),h(0))$. Then $\beta<1$ since as we saw in the last proof, $h(0)\leq 1$. 
We also saw in that proof that $g_2''(z_0)\neq 0$, and since $g_1''(x)<0$ when $x<z_0$,
$h(x)\rightarrow -\infty$ as $x\uparrow z_0$. Therefore,
there is a unique $z_\beta<z_0$ so that $h(z_\beta)=\beta$.
By Lemma A.1(g), as a function of $x$, $\gamma(\beta,x)$  increases on $(0,z_\beta)$
and decreases on $(z_\beta,x_1^*)$. By Lemma A.2(a), the maximum
$\gamma(h(z_\beta),z_\beta)<0$. But by Lemma A.1(f) and (\ref{21}),
$\gamma(\beta,\bar{u}^*)>\gamma(1,\bar{u}^*)=0$, a contradiction since
$\bar{u}^*<x_1^*<\underline{z}$.

For the last part, assume that $\beta\in (h(\underline{z}),1]$. Assume first
that $x\in (z_0,h_-^{-1}(\beta))$. Then $g_1''(x)>0$ and
$h(x)>h(h_-^{-1}(\beta))=\beta$. Therefore,
$$g''(x;\beta)=g_1''(x)(h(x)-\beta)>0.$$
Assume next that $z_0>0$ and that $x\in (0,z_0)$. Then $g_1''(x)<0$ and
$h(x)<h(0)\leq h(\underline{z})<\beta$, where we used the result just proved.
Therefore,  $g''(x;\beta)>0$ again. Finally, assume that $x=z_0$. Then
$g''(z_0;\beta)=g_2''(z_0)=g''(z_0)>0$ since $z_0<x_1^*$ and $g_1''(z_0)=0$.
This ends the proof.

\vspace{0.2cm}
\b{Proof of Lemma 2.4.} For any $\beta\in (h(\underline{z}),1)$, by Lemma A.1(f) and (\ref{21}), 
$\gamma(\beta,\bar{u}^*)>\gamma(1,\bar{u}^*)=0$, and by Lemma A.2(a), $\gamma(\beta,h_-^{-1}(\beta))<0$.
Therefore, since $\gamma(\beta,x)$ is continuous in $x$, there is at least one 
$u_1(\beta)\in (\bar{u}^*,h_-^{-1}(\beta))$ so that 
$\gamma(\beta,u_1(\beta))=0$. Since $u_1(\beta)<h_-^{-1}(\beta) $,
it follows from Lemma A.3 that $g''(u_1(\beta);\beta)>0$, and so by Lemma A.2(g),
$$\frac{\partial}{\partial x}\gamma(\beta,u_1(\beta))<0$$
for all $\beta\in (h(\underline{z}),1)$. This implies in particular that $u_1(\beta)$ is unique,
and furthermore by the implicit function theorem, $u_1(\beta)$ is continuously differentiable on
$(h(\underline{z}),1)$. Finally, using that $\frac{\partial}{\partial \beta}\gamma(\beta,u_1(\beta))=0$
and Lemma A.2(g) yields
$$u_1'(\beta)=-\frac{\frac{\partial}{\partial \beta}\gamma(\beta,u_1(\beta))}
{\frac{\partial}{\partial x}\gamma(\beta,u_1(\beta))}<0,$$
and so $u_1$ is strongly decreasing.

\vspace{0.2cm}
\b{Proof of Lemma 2.5.} Taking the derivative w.r.t $\beta$  in $\gamma(\beta,u_1(\beta))=0$ gives after some
simplification that
\begin{equation}\frac{d}{d \beta}g'(u_1(\beta);\beta)=-\frac{g_1(u_1(\beta))g'(u_1(\beta);\beta)}
{g(u_1(\beta);\beta)}.\label{a4}\end{equation}
Taking the derivative of $G(\beta)$ and using (\ref{a4}) gives
\begin{equation}
G'(\beta)=\frac{g_{2,\infty}'g_1(u_1(\beta))-g_{1,\infty}'g_2(u_1(\beta))}
{g'(u_1(\beta);\beta)g(u_1(\beta);\beta)}.\label{a5}
\end{equation}
We have
$$\frac{g_2'(x)}{g_1'(x)}-\frac{g_2(x)}{g_1(x)}=\frac{W(x)}{g_1'(x)g_1(x)},$$
which is positive for $x>x_\lambda$. Consequently, by Lemma A.1(b),
$$\frac{g_2'(x)}{g_1'(x)}>\frac{g_2(x)}{g_1(x)}>\frac{g_2(u_1(\beta))}{g_1(u_1(\beta))},
\;\;\;\;x>x_\lambda.$$
Letting $x\rightarrow\infty$ gives
$$\frac{g_{2,\infty}'}{g_{1,\infty}'}>\frac{g_2(u_1(\beta))}{g_1(u_1(\beta))},$$
and the result follows from (\ref{a5}).

\vspace{0.2cm}
\b{Proof of Lemma 2.6.} The part about $u_1$ is proved in Lemma 2.4.
Since $g'(x;\beta)$ has a local maximum at $x=h_-^{-1}(\beta)$ and a local minimum at
$x=h_+^{-1}(\beta)$, for existence of $u_2$ and $u_3$ all we need to show is that for
all $\beta\in (\alpha_0,\beta_0)$,
\begin{enumerate}
\item $g'(u_1(\beta);\beta)>g'(h_+^{-1}(\beta);\beta).$
\item $\lim_{x\rightarrow\infty}g'(x;\beta)>g'(u_1(\beta);\beta).$
\end{enumerate}
But (1) is satisfied by the definition of $\beta_0$ in (\ref{28}) since $\beta<\beta_0$, while (2)
follows since $G(\alpha_0)\geq 1$ and $G$, whenever finite,  is strictly increasing by Lemma 2.5.
Furthermore, letting $w(\beta,u)=g'(u_1(\beta);\beta)-g'(u;\beta)$ gives
$w(\beta,u_2(\beta))= w(\beta,u_3(\beta))=0$, and then it easily follows from the implicit 
function theorem that $u_2$ and $u_3$ are continuously differentiable.

When $\alpha_0=h(\underline{z})$ it follows from Lemma 2.4 that $u_1(\alpha_0)=\underline{z}=h_-^{-1}(\alpha_0)$
and so $g'(x;\alpha_0)$ has a maximum at $x=\underline{z}$, which implies that
$u_1(\alpha_0)=u_2(\alpha_0)=\underline{z}$.

When $\alpha_0>h(\underline{z})$, then since $G(\alpha_0)=1$,
$$\frac{\lim_{x\rightarrow\infty}g'(x;\beta)}{g'(u_1(\beta);\beta)}\downarrow  1
\;\;\;\mbox{as}\;\;\beta\downarrow \alpha_0.$$
Therefore, to find the $u_3(\beta)$ that satisfies $g'(u_3(\beta);\beta)=
g'(u_1(\beta);\beta)$ we have to go further and further out as $\beta\downarrow \alpha_0$.

If $\beta_0<\min\{1,h_\infty\}$, by continuity $v(\beta_0)=1$, and so we can set $u_3(\beta_0)=u_2(\beta_0)=h_+^{-1}(\beta_0).$

\vspace{0.2cm}

\b{Proof of Lemma 2.7.} Using that $\gamma(\beta,u_1(\beta))=0$ 
and that $g'(u_2(\beta);\beta)=g'(u_1(\beta);\beta)$
we get
$$J_1(\beta)=J_1(\beta)- \gamma(\beta,u_1(\beta))=  k\left(u_2(\beta)-\frac{g(u_2(\beta);\beta)}
{g'(u_1(\beta);\beta)}\right)-K.$$

Since
$$\frac{d}{d \beta}g(u_2(\beta);\beta)=g'(u_1(\beta);\beta)u_2'(\beta)-g_1(u_2(\beta))$$
we get
$$\frac{d}{d \beta}J_1(\beta)= \frac{k}{(g'(u_1(\beta);\beta))^2}v_1(\beta),$$
where
\begin{eqnarray*}
v_1(\beta)&=&g'(u_1(\beta);\beta)g_1(u_2(\beta))+g(u_2(\beta);\beta)\frac{d}{d \beta}
g'(u_1(\beta);\beta)\\
&=&g'(u_1(\beta);\beta)g_1(u_2(\beta))-\frac{g(u_2(\beta);\beta)g_1(u_1(\beta))g'(u_1(\beta);\beta)}
{g(u_1(\beta);\beta)},
\end{eqnarray*}
where we used (\ref{a4}) in the second equality. Clearly $J_1(\beta)$ and $v_1(\beta)$ have
the same sign, and 
since $g(u_1(\beta);\beta)$ and $g'(u_1(\beta);\beta)$ are both positive,
the sign of $v_1(\beta)$ is the same as the sign of
\begin{eqnarray*}
v_2(\beta)&=&g_1(u_2(\beta))g(u_1(\beta);\beta)-g_1(u_1(\beta))g(u_2(\beta);\beta)\\
&=&g_1(u_2(\beta))g_2(u_1(\beta))-g_2(u_2(\beta))g_1(u_1(\beta))\\
&=&g_1(u_1(\beta))g_1(u_2(\beta))\left(\frac{g_2(u_1(\beta))}{g_1(u_1(\beta))}-
\frac{g_2(u_2(\beta))}{g_1(u_2(\beta))}\right)<0,
\end{eqnarray*}
by Lemma A.1(b). Replacing $u_2(\beta)$ by $u_3(\beta)$ we get the same conclusion for
$J_{13}(\beta)$.

\vspace{0.2cm}
\b{Proof of Proposition 2.1.} Since either R5 or R6 apply, it follows from
the comment after Lemma 2.3 that $h(\underline{z})<h_\infty$. Therefore, by
Lemma 2.6, for any $\beta\in [\alpha_0,\beta_0)$ there are numbers
$u_i=u_i(\beta)\,i=1,2,3$, so that
$$g'(u_1;\beta)=g'(u_2;\beta)=g'(u_3;\beta),$$
where $u_1=u_2=\underline{z}$ if $\beta=\alpha_0=h(\underline{z})$.

Assume G1 or G2. We will show that
\begin{equation}
J_{13}(\alpha_0)\geq 0\;\;\;\mbox{and}\;\;\;J_{13}(\beta_0)<0.
\label{a51}\end{equation}
In case G1, $\alpha_0=h(\underline{z})$ and by (\ref{212}),
$J_1(h(\underline{z}))=0$. Also
$J_2(h(\underline{z}))>0$ and so $J_{13}(\alpha_0)=J_1(\alpha_0)+
J_2(\alpha_0)>0$. In case G2, $J_{13}(\alpha_0)\geq 0$ by assumption.
To prove that $J_{13}(\beta_0)<0$ assume first that $\beta_0<1$, which
implies that $\beta_0<\min\{1,h_\infty\}$ by Lemma 2.2. Therefore,
$J_1(\beta_0)<0$ while $J_2(\beta_0)=0$ by (\ref{214}), and so
$J_{13}(\beta_0)<0$. If $\beta_0=1$, then  $J_{13}(1)=K_3-K_2< 0$,
and so (\ref{a51}) is proved.

By Lemma 2.7, $J_{13}$ is strongly decreasing on $(\alpha_0,\beta_0)$,
and so by (\ref{a51}) there is a unique $\hat{\beta}$ so that 
$J_{13}(\hat{\beta})=0$, i.e. so that
$-J_1(\hat{\beta})=J_2(\hat{\beta})$. Assume first that $\beta_0<1$.
If $J_2(\hat{\beta})\leq K$,
let $\tilde{\beta}=\hat{\beta}$. If $J_2(\hat{\beta})>K$, 
since $J_2(\beta_0)=0$, we can define $\tilde{\beta}<\beta_0$ by
$$\tilde{\beta}=\min\{\beta>\hat{\beta}:J_2(\beta)=K\}.$$
Then  since $J_{13}(\beta)$ is decreasing,
$$0>J_{13}(\tilde{\beta})=K+J_1(\tilde{\beta}),$$
hence $-J_1(\tilde{\beta})>K.$

If $\beta_0=1$, then  $J_2(1)=K_3<K$ and we can use the same arguments. 

Finally assume that G3 holds. Then $J_{13}(\alpha_0)<0$, but since $J_2(\alpha_0)>K$ and  
$$J_2(\beta_0)=\left\{\begin{array}{lll}
0,\;\;\;&\beta_0<1,\\
K_3,&\beta_0=1,
\end{array}\right.$$
and $K_3<K$, we can define $\tilde{\beta}<\beta_0$ by
$$\tilde{\beta}=\min\{\beta>\alpha_0:J_2(\beta)=K\}.$$
Since $0>J_{13}(\alpha_0)>J_{13}(\tilde{\beta})$, $-J_1(\tilde{\beta})>J_2(\tilde{\beta})$, and so
case (iii-b) applies.

\vspace{0.2cm}
\b{Proof of Proposition 2.2.} Assume first that $h_\infty=1$. 
By letting $h_+^{-1}(1)=
\lim_{\beta\rightarrow 1}h_+^{-1}(\beta)   =\infty$
and using that $u_1(\beta)\rightarrow \bar{u}^*$ as $\beta\rightarrow 1$ we can use (\ref{27}) to
conclude that $v(\beta)\rightarrow v(1)$  as $\beta\rightarrow 1$, where
$$v(1)=\frac{g_{\infty}'}{g'(\bar{u}^*)}=\frac{c^*}{c^{\infty}}>1,$$
by assumption.
Therefore,  $\beta_0<1$, and
so trivially $\beta_0<\min\{1,h_\infty\}$. Clearly $\alpha_0$ can be well defined as well,
and $\alpha_0<\beta_0$ by (\ref{29}). Hence the assumptions of Lemma 2.6 hold.

Assume next  that $h(\underline{z})<h_\infty<1$ and that $G(h_\infty)>1$. We have by (\ref{29}),
$$1<\frac{G(\beta)}{v(\beta)}=\frac{\lim_{x\rightarrow\infty}g'(x;\beta)}
{g'(h_+^{-1}(\beta);\beta)}\rightarrow 1\;\;\;\mbox{as}\;\;\beta\uparrow h_\infty,$$
and since $G(h_\infty)>1$, it follows that $v(\beta)>1$ for $\beta$ sufficiently close to
$h_\infty$. Therefore,
$\beta_0<h_\infty=\min\{1,h_\infty\}$ and
again all assumptions of Lemma 2.6 hold.

The rest of the proof is now the same as the proof of Proposition 2.1
for the case with $\beta_0<1$. 

\vspace{0.2cm}
We will now turn to the proof of Theorem 3.1. The proof is standard once the necessary 
variational inequalities have been established, and that is the topic of the next two lemmas.

\begin{lemma}
Let $V$ be as the proposed value functions in Theorem 3.1. Then $V$ is continuously differentiable on
$(0,\infty)$ and twice continuously differentiable on the set $(0,\bar{u}^*)\cup 
(\bar{u}^*,\tilde{u}_1)\cup (\tilde{u}_1,\tilde{u}_3)\cup (\tilde{u}_3,\infty)$.
Furthermore, 
$$LV(x)= 0\;\;\mbox{on}\;\;(0,\bar{u}^*)\cup(\tilde{u}_1,\tilde{u}_3)\;\;\;
\mbox{and}\;\;\;LV(x)<0\;\;\mbox{on}\;\;(\bar{u}^*,\tilde{u}_1)\cup(\tilde{u}_3,\infty).$$

\end{lemma}
{\sl Proof.}  Consider first case (iii-a). By definition, 
$c^*g(\bar{u}^*-)=k\bar{u}^*-K$ and $c^*g'(\bar{u}^*-)=k$, hence $V$ is continuously differentiable at
$\bar{u}^*$. Also, by definition
$$V(\tilde{u}_1)=\frac{k}{g'(\tilde{u}_1;\tilde{\beta})}
g(\tilde{u}_1;\tilde{\beta})
=\gamma(\tilde{\beta},\tilde{u}_1)
+k\tilde{u}_1-K=k\tilde{u}_1-K=V(\tilde{u}_1-)  $$
and clearly $V'(\tilde{u}_1+)=k$, hence $V$ is continuously differentiable at $\tilde{u}_1$.
Similarly, $V$ is continuously differentiable at $\tilde{u}_3$.

By definition, $LV(x)=0$ for $x\in (0,\bar{u}^*)\cup(\tilde{u}_1,\tilde{u}_3)$. Let
$x\in (\bar{u}^*,\tilde{u}_1)$. Then
$$LV(x)=\mu(x)k-\lambda(kx-K)<0,$$
by (\ref{25}) and the fact that $\tilde{u}_1<\underline{z}$.

Finally, let $x\in (\tilde{u}_3,\infty)$. Then
\begin{eqnarray*}
LV(x)&=&\mu(x)V'(x)-\lambda V(x)\\
&\leq  & \mu(\tilde{u}_3)k-\lambda V(\tilde{u}_3)\\
&=&\mu(\tilde{u}_3-)V'(\tilde{u}_3-)-\lambda V(\tilde{u}_3-)\\
&< & \frac 12 \sigma^2(\tilde{u}_3-)V''(\tilde{u}_3-)+
\mu(\tilde{u}_3-)V'(\tilde{u}_3-)-\lambda V(\tilde{u}_3-)\\
&=&LV(\tilde{u}_3-)=0.
\end{eqnarray*}
Here the first inequality follows from the fact that for $x>\tilde{\mu}_3$,
$\frac{d}{dx}(\mu(x)k-\lambda V(x))=k(\mu'(x)-\lambda)\leq 0$. The next equality is clear
since $V$ is continuously differentiable and $V'(\tilde{u}_3+)=k$. For the second inequality,
note that as in (\ref{26}),
$$g''(\tilde{u}_3;\tilde{\beta})=g_1''(\tilde{u}_3)(h(\tilde{u}_3)-\tilde{\beta}).$$
But $g_1$ is convex on $(x_{\lambda},\infty)$ and $\tilde{\beta}<h(\tilde{u}_3)$ by Proposition 2.1
since $h_+$ is increasing.
Consequently $V''(\tilde{u}_3-)>0$.

Now to case (iii-b). That $V$ is continously differentiable at $\bar{u}^*$ and $\tilde{u}_1$ is
proved as above.  Furthermore,
\begin{eqnarray*}
V(\tilde{u}_3)&=&c_{\tilde{\beta}}g(\tilde{u}_2;\tilde{\beta})+k(\tilde{u}_3-\tilde{u}_2)-K\\
&=&c_{\tilde{\beta}}g(\tilde{u}_3;\tilde{\beta})+J_2(\tilde{\beta})-K\\
&=&c_{\tilde{\beta}}g(\tilde{u}_3;\tilde{\beta})\\
&=&V(\tilde{u}_3-).
\end{eqnarray*}
Trivially, $V'(\tilde{u}_3+)=V'(\tilde{u}_3-)=k$.
Finally, the signs of $LV(x)$ are shown  just as above. 

\vspace{0.3cm}
For a function $\phi:[0,\infty)\mapsto
[0,\infty),$  define the maximum utility operator $M$ by
\begin{eqnarray*}
M\phi(x)&=&\!{\rm sup}\{\phi(x-\eta)+k\eta-K: 0\leq \eta \leq x\}\\
&=&\!{\rm sup}\{\phi(z)+k(x-z)-K: 0\leq z \leq x\}.
\end{eqnarray*}

\begin{lemma}
Let $V$ be as in Lemma A.3. Then 
$$MV(x)=V(x)\;\;\mbox{on}\;\;[\bar{u}^*,\tilde{u}_1]\cup [\tilde{u}_3,\infty)
\;\;\;\mbox{and}\;\;\;MV(x)<V(x)\;\;\mbox{on}\;\;(0,\bar{u}^*)\cup(\tilde{u}_1,\tilde{u}_3).$$

\end{lemma}
{\sl Proof.} 
  The following elementary argument will often be used in the proof.
Assume we want to maximize $h(y)=\phi(y)+k(x-y) -K$ for $y\in [a_0,a_1]$ with $0\leq a_0<a_1\leq x$.
Then $h'(y)=\phi'(y)-k$, so if $\phi'(y)\leq k$ on $[a_0,a_1]$,
a maximum is at $y=a_0$. Conversely, if $\phi'(y)\geq k$ on $[a_0,a_1]$,
a maximum is at $y=a_1$.

Consider first case (iii-a).
For $x\in (0,\tilde{u}_1]$,  $V'(y)\leq k$ for $y\in (0,x)$, hence $MV(x)=V(0)+kx-K=kx-K$. Therefore,
$MV(x)=V(x)$ on $[\bar{u}^*,\tilde{u}_1]$, while for $x\in [0,\bar{u}^*)$,
$$MV(x)=(V(\bar{u}^*)-k\bar{u}^*+K)+kx-K= V(x)+\int_x^{\bar{u}^*}(V'(y)-k)dy<V(x).$$
Let $x\in (\tilde{u}_1,\tilde{u}_2)$. Then 
$$V'(x)\left\{\begin{array}{lll}<0,\;\;\;&y\in (0,\tilde{u}_1),\\
>0, & y\in (\tilde{u}_1,x).\end{array}\right.$$
Therefore, $MV(x)=\max\{V(x)-K,kx-K\}$. Clearly $V(x)-K<V(x)$ and
\begin{eqnarray*}
V(x)-(kx-K)=\int_{\tilde{u}_1}^x(V'(y)-k)dy+V(\tilde{u}_1)-(k\tilde{u}_1-K)
=\int_{\tilde{u}_1}^x(V'(y)-k)dy>0,
\end{eqnarray*}
and so $MV(x)<V(x)$.

Let $x\in [\tilde{u}_2,\tilde{u}_3)$. Then, since $V'(y)>k$ on $(\tilde{u}_1,\tilde{u}_2)$ 
while $V'(y)< k$ on $(0,\tilde{u}_1)\cup (\tilde{u}_2,x)$, we get
\begin{equation}MV(x)=\max\{V(\tilde{u}_2)+k(x-\tilde{u}_2)-K,kx-K\}.
\label{a6}\end{equation}
Now,
$$V(\tilde{u}_2)-k\tilde{u}_2=V(\tilde{u}_1)-k\tilde{u}_1-
J_1(\tilde{\beta})\leq V(\tilde{u}_1)-k\tilde{u}_1+K=0,$$
and therefore $V(\tilde{u}_2)+k(x-\tilde{u}_2)-K\leq kx-K$, implying that   
$MV(x)=kx-K$. But then
$$MV(x)-V(x)=k(x-\tilde{u}_1)-(V(x)-V(\tilde{u}_1))=J(\tilde{\beta},\tilde{u}_1,x)
< J(\tilde{\beta},\tilde{u}_1,\tilde{u}_3)=0,$$
and consequently $MV(x)<V(x)$.

Finally, let $x\in [\tilde{u}_3,\infty)$. Then $MV(x)$ is again given by
(\ref{a6}), and as there $MV(x)=kx-K=V(x)$.

We continue with case (iii-b). When $x\in [0,\tilde{u}_2)$ the proof is just as above,
so assume that $x\in [\tilde{u}_2,\tilde{u}_3)$. Then again $MV(x)$ is given by
(\ref{a6}), but now
$$V(\tilde{u}_2)-k\tilde{u}_2=V(\tilde{u}_1)-k\tilde{u}_1-
J_1(\tilde{\beta})> V(\tilde{u}_1)-k\tilde{u}_1+K=0,$$
so
$MV(x)=V(\tilde{u}_2)+k(x-\tilde{u}_2)-K$ and therefore
$$MV(x)-V(x)=J(\tilde{\beta},\tilde{u}_2,x)-K<J(\tilde{\beta},\tilde{u}_2 ,\tilde{u}_3)-K=0,$$
and so $MV(x)<V(x)$.

When $x\in [\tilde{u}_3,\infty)$, as above $MV(x)=V(\tilde{u}_2)+k(x-\tilde{u}_2)-K=V(x).$

\vspace{0.2cm}
\b{Proof of Theorem 3.1.} Using the variational inequalities in Lemmas A.3 and A.4, the proof is exactly as
the proof of Theorem 2.2 in \cite{Bai10b}.

\vspace{0.2cm}
\b{Proof of Theorem 3.2.} The proof is exactly like the proof of the (more interesting) Theorem 4.2 below.

\vspace{0.2cm}
\b{Proof of Theorem 3.3.} We can set $u_3(1)=\infty$ and then since it is case R6, $J_{13}(1)=K_3-K_2<0$ and
$J_2(1)=K_3<\min\{K_2,K\}$. With $V(x)=c^*g(x)$, clearly $LV(x)=0$. The proof of the
variational inequality $MV(x)<V(x)$ is  as in
Lemma A.4. Then the rest is as in the proof of Theorem 4.3 below.

\vspace{0.2cm}
\b{Proof of Theorem 4.1.} That is exactly like the proof of Theorem 3.1.

\vspace{0.2cm}
\b{Proof of Lemma 4.2.} Since $G(1)=\frac{g_{\infty}'}{g'(\bar{u}^*)}=
\frac{c^*}{c^\infty}> 1$ by assumption,
it follows from Lemma 2.5 that if it exists, $\hat{\beta}< 1$.
In case (1), note that by Lemma 2.4, $g'(u_1(h(\underline{z}));h(\underline{z}))=
g'(\underline{z};h(\underline{z}))$ and by (\ref{26}),
$$g''(x;h(\underline{z}))=g_1''(x)(h(x)-h(\underline{z}))<0,\;\;\;x>\underline{z},$$
since $\underline{z}>z_0$ and $h(\underline{z})\geq h_\infty$. Therefore, 
$$\lim_{x\rightarrow\infty}g'(x;h(\underline{z}))<g'(u_1(h(\underline{z}));h(\underline{z})),$$
i.e. $G(h(\underline{z}))<1$ so $\hat{\beta}$ exists 
and $\hat{\beta}>h(\underline{z})$. In case (2) $G(h_\infty)\leq 1$ by assumption, and the result follows.
Also $\hat{\beta}\geq h_\infty >h(\underline{z})$.
In case (3),  $G(\alpha_0)=1$, hence $\hat{\beta}=\alpha_0$.  Since in all three cases 
$\hat{\beta}>h(\underline{z})$,
$$u_1(\hat{\beta})<h_-^{-1}(\hat{\beta})<h_-^{-1}(h(\underline{z}))=\underline{z},$$ 
and this ends the proof.

\vspace{0.2cm}
\b{Proof of Theorem 4.2.} 
The proof that the proposed value function $V$ (we write $V$ instead ov $V^*$) is 
continuously differentiable is straightforward, and is omitted.

We must prove the variational inequalities
\begin{equation}
LV(x)=0\;\;\mbox{on}\;\;(0,\bar{u}^*)\cup (\hat{u}_1,\infty)
\;\;\;\mbox{and}\;\;\;LV(x)<0\;\;\mbox{on}\;\;(\bar{u}^*,\hat{u}_1).\label{a8}
\end{equation}
The first part is obvious, and the second is identical to that in the proof of Lemma A.3.

The next step is to prove the variational inequalities
\begin{equation}
MV(x)<0\;\;\mbox{on}\;\;(0,\bar{u}^*)\cup (\hat{u}_1,\infty)
\;\;\;\mbox{and}\;\;\;MV(x)=0\;\;\mbox{on}\;\;(\bar{u}^*,\hat{u}_1).\label{a9}
\end{equation}

We begin by proving
\begin{equation}
c_{\hat{\beta}}g'(x;\hat{\beta})>k,\;\;\,x>\hat{u}_1,\label{a7}
\end{equation}
where we throughout the proof write $\hat{u}_1$ for $u_1(\hat{\beta})$.
This is equivalent to
$$m(x)\defi g'(x;\hat{\beta})-g'(\hat{u}_1;\hat{\beta})>0,\;\;\;x>\hat{u}_1.$$
Note that $m(\hat{u}_1)=0$ and since $G(\hat{\beta})=1$, $\lim_{x\rightarrow\infty}m(x)=0$.
However, $m'(x)=g''(x;\hat{\beta})=g_1''(x)(h(x)-\hat{\beta})$, so in particular
$$m'(\hat{u}_1)=g_1''(\hat{u}_1)(h(\hat{u}_1)-\hat{\beta}).$$
Clearly $g_1''(\hat{u}_1)>0$. Also, $\hat{u}_1<h_-^{-1}(\hat{\beta})$,
and since $h$ is decreasing on $(x_1^*,x_\lambda)$, $h(\hat{u}_1)>
h(h_-^{-1}(\hat{\beta}))=\hat{\beta}$. Consequently, $m'(\hat{u}_1)>0$. Furthermore,
$m'(x)$ has only one zero on $(\hat{u}_1,\infty)$ and so (\ref{a7}) is proved.
To prove (\ref{a9}), we only prove that $MV(x)<0\;\;\mbox{on}\;\; (\hat{u}_1,\infty)$, since the rest is as in Lemma A.4. We get for $x>\hat{u}_1$,
$$MV(x)={\rm max}\{V(x)-K, kx-K\}.$$
Now
\begin{eqnarray*}
kx-K-V(x)&=&\int_{0}^{x}(k-V'(y))dy-K\\
&=&\int_{0}^{\bar{u}^{*}}(k-c^*g'(y))dy+\int_{\hat{u}_1}^{x}(k-c_{\hat{\beta}}g'(y;\hat{\beta}))dy-K\\
&=&\int_{\hat{u}_1}^{x}(k-c_{\hat{\beta}}g'(y;\hat{\beta}))dy<0,
\end{eqnarray*}
where we used (\ref{a7}) in the last inequality.
Therefore, $MV(x)<V(x)$ on $(\hat{u}_1,\infty)$.

Using the variational inequalities (\ref{a8}) and (\ref{a9}), standard arguments, see e.g. \cite{Bai10b}, shows that
$V_{\bar{u}^{*},0}(x)=V(x)=V^*(x)$ for $x\leq \hat{u}_1$.

For the lump sum strategy dividend barrier strategy $\pi_{\bar{u},0}$ with $\bar{u}>\hat{u}_1$, it is easy to prove that its value function is
\begin{eqnarray*}
V_{\bar{u},0}(x)=\left\{\begin{array}{llll}
\delta(\bar{u})g(x),\;\;\;&x\in[0,\bar{u}),\\
kx-K,& x\in [\bar{u},\infty),\\
\end{array}\right.
\end{eqnarray*}
where $\delta(\bar{u})=\frac{k\bar{u}-K}{g(\bar{u})}$. Then, by the same  arguments as in Theorem 2.1(c) in \cite{Bai10b},  $V_{\bar{u},0}(x)$ is increasing and
${\mathop{\lim}_{\bar{u}\rightarrow
\infty}V_{\bar{u},0}(x)}=V(x)$ for $x>\hat{u}_1$.\\
 \indent  Assume that there exists an optimal strategy $\pi^{*}$
 when the initial surplus $x>\hat{u}_1$.   
Since $V$ is twice continuously differentiable except at a finite number of points,
It\^{o}'s formula used on 
$e^{-\lambda t}V(X_{t+})$ stopped at the time of ruin $\tau^{\pi^*}$ gives
\begin{eqnarray*}
e^{-\lambda (t\wedge \tau^{\pi^*})}V(X_{(t\wedge \tau^{\pi^*})+}^{\pi^*})
&=& V(x)+\int_0^{t\wedge \tau^{\pi^*}}
e^{-\lambda s}LV(X_s^{\pi^*})ds \nonumber\\
&&+\int_0^{t\wedge \tau^{\pi^*}}
e^{-\lambda s}\sigma(X_s^{\pi^*})V'(X_s^{\pi^*})dW_s \nonumber \\
&&+\sum_{s\leq t\wedge \tau^{\pi^*}} e^{-\lambda s} \left(V(X_{s+}^{\pi^*})-V(X_{s}^{\pi^*})
\right) \\
&<&  V(x)+\int_0^{t\wedge \tau^{\pi^*}}
e^{-\lambda s}\sigma(X_s^{\pi^*})V'(X_s^{\pi^*})dW_s \nonumber \\
&&-\sum_{n=1}^{\infty}e^{-\lambda \tau_n^{\pi^*}}(k\xi_n^{\pi^*}-K)1_{\{
\tau_n^{\pi^*}\leq t\wedge \tau^{\pi^*}\}}, \nonumber
\end{eqnarray*}
The inequality follows from the fact that if $X^{\pi^{*}}_{\tau_{1}^{*}}<\hat{u}_1$
then by (\ref{a8}),
$$P_x\left(\int_0^{t\wedge \tau^{\pi^*}}
e^{-\lambda s}LV(X_s^{\pi^*})ds<0\right)=1,$$
while if $X^{\pi^{*}}_{\tau_{1}^{*}}\geq \hat{u}_1$ then by (\ref{a9}),
$$P_x\left(\sum_{s\leq t\wedge \tau^{\pi^*}} e^{-\lambda s} \left(V(X_{s+}^{\pi^*})-V(X_{s}^{\pi^*})\right)
<-\sum_{n=1}^{\infty}e^{-\lambda \tau_n^{\pi^*}}(k\xi_n^{\pi^*}-K)
1_{\{\tau_n^{\pi^*}\leq t\wedge \tau^{\pi^*}\}}\right)=1.$$
Taking expectations and letting $t\rightarrow\infty$ gives
$0<V(x)-V_{\pi^*}(x)$, a contradiction since ${\mathop{\lim}_{\bar{u}\rightarrow
\infty}V_{\bar{u},0}(x)}=V(x)$ for $x>\hat{u}_1$.

\vspace{0.2cm}
\b{Proof of Theorem 4.3.} The proof is basically the same as that of Theorem 3.3, and is omitted.


\begin{thebibliography}{9}

\bibitem{Al98}
\sc{L.~H. R. Alvarez} 
{\em Optimal harvesting under stochastic fluctuations and critical depensation},
Math. Biosci., 152, (1998), pp.~63--85.

\bibitem{Bai10b}
\sc{L.~Bai and J.~Paulsen.} {\em Optimal dividend policies with
transaction costs for a class of diffusion processes.} 
SIAM J. Control Optim., 48, (2010), pp.~4987--5008.

\bibitem{Bai11}
\sc{L.~Bai, M.~Hunting and J.~Paulsen.} {\em Optimal 
dividend policies for a class of growth-restricted diffusion
processes  under
transaction costs and solvency constraints.}
To appear in Finance Stoch.


\bibitem{Bor} 
\sc{A.~Borodin and P.~Salminen.} {\em Handbook of Brownian Motion: Facts and Formulae, 2nd rev. edn. Probability
and Its Applications.} Birkh\"{a}user, Basel, 2002.

\bibitem{Eg08}
\sc{M.~Egami.} {\em A direct solution method for stochastic
impulse control problems of one-dimensional diffusions.}
SIAM J. Control Optim., 47, (2008), pp.~1191--1218.

\bibitem{P3}
\sc{J.~Paulsen.}  {\em Optimal dividend payments until ruin of 
diffusion processes when payments are subject to both fixed and 
proportional costs.} Adv. Appl. Prob., 39, (2007), pp.~669--689.

\bibitem{P4}
\sc{J.~Paulsen.}   {\em Optimal dividend payments and reinvestments of 
diffusion processes when payments are subject to both fixed and 
proportional costs.} SIAM J. Control Optim., 47 (2008), pp.~2201--2226.


\end{thebibliography}
\end{document}